\documentclass[a4paper]{article}

\usepackage{amsfonts}
\usepackage{amssymb, amscd}
\usepackage{latexsym}

\newtheorem{theorem}{Theorem}[section]

\newtheorem{prop}[theorem]{Proposition}

\newtheorem{prob}[theorem]{Problem}
\newtheorem{remark}[theorem]{Remark}
\newtheorem{definition}[theorem]{Definition}
\newtheorem{example}[theorem]{Example}
\newtheorem{problem}[theorem]{Problem}
\newtheorem{defin}[theorem]{Definition}

\def\ol{\overline}

\def\mf{\mathfrak}
\def\kp{Val}

\begin{document}


\title{Algebraic logic and logically-geometric types in varieties of algebras }

\maketitle

\begin{center}

\author{B. Plotkin$^{a}$,
     \and E. Aladova$^{b}$,
      \and E. Plotkin$^{c}$}

 \smallskip
        {\small

        $^{a}$ Department of Mathematics,
          Hebrew University of Jerusalem,
          91904, Jerusalem, Israel

                {\it E-mail address:} plotkin@macs.biu.ac.il

        }

 \smallskip
        {\small
               $^{b}$ Department of Mathematics,
          Bar-Ilan University,

          52900, Ramat Gan, Israel

                {\it E-mail address:} aladovael@mail.ru
        }

  \smallskip
        {\small
               $^{c}$ Department of Mathematics,
          Bar-Ilan University,

          52900, Ramat Gan, Israel

                {\it E-mail address:} plotkin@macs.biu.ac.il
        }

\end{center}

\begin{abstract}
The main objective of this paper is to show that the notion of
type which was developed within the frames of logic and model
theory has deep ties with geometric properties of algebras.  These
ties go back and forth from universal algebraic geometry to the
model theory through the machinery of algebraic logic. We show
that types  appear naturally as logical kernels in the Galois
correspondence between filters in the Halmos algebra of first
order formulas with equalities and elementary sets in the
corresponding affine space.
\end{abstract}


\tableofcontents

\section{Introduction}\label{S_Int}

The main objective of the paper is to show that the notion of
type which was developed within the frames of logic and model
theory has deep ties with geometric properties of algebras.  These
ties go back and forth from universal algebraic geometry to model
theory through the machinery of algebraic logic.

More precisely, we shall show that types appear naturally as
logical kernels in the Galois correspondence between filters in
the Halmos algebra of first order formulas with equalities and
elementary sets in the corresponding affine space. Note that in
our terminology  the term "elementary set" has the meaning of
"definable set" in the standard model theoretic terminology.  This
Galois correspondence generalizes classical Galois correspondence
between ideals in the polynomial algebra and algebraic sets in the
affine space. The sketch of the ideas of universal algebraic
geometry can be found in \cite{Plotkin_AG}, \cite{Plotkin_7Lec},
\cite{Plotkin_SomeResultsUAG}, \cite{PlotZitom_1},
\cite{PlotZitom_2}, \cite{BMR}, \cite{MR},
 \cite{KvasMyasnRemes_3}, \cite{DanMyasnRem_1} \cite{DanMyasnRem_2},
\cite{DanMyasnRem_3}, \cite{DanKazachRem_1},
\cite{DanKazachRem_2}, \cite{KhM_1}, \cite{KhM_2},
\cite{RipsSela}, \cite{Sela1},   etc. As for standard definitions
of model theory, we refer to monographs \cite{Marker},
\cite{Poizat}, \cite{ChangKeisler_ModelTh}, \cite{Hodges},  etc.
For the exposition of concepts and  results of algebraic logic see
\cite{Halmos} -- \cite{Halmos_4},  \cite{HMT},
\cite{FontJanPig_1}, \cite{FontJanPig_2}, \cite{AndNemSain},
\cite{AmerAhmed}, etc.

Methodologically, in the paper we give a sketch of some  ideas
which provide interactions of algebraic logic with geometry, model
theory and algebra. We believe that a  development of the
described approach can make benefits to each of these areas. We
shall stress that the paper does not contain a bunch of new
results. Its main duty is to specialize  new problems and  to
underline  common points of algebra, logic and geometry through
the notion of the type.

The paper is organized as follows. Section \ref{sec:al} is devoted
to structures of algebraic logic.
 We define here various kinds  of Halmos algebras, consider the value homomorphism
and provide the reader with the main examples of algebras under consideration.  Section \ref{sec:galois} deals with basic approaches of universal algebraic geometry. We define the general Galois correspondence which plays the important role in all considerations. The description of this correspondence starts from the classical case and extends to the case of multi-sorted logical geometry over an arbitrary variety of algebras. In  Section \ref{sec:mod} we  recall the model theoretic notion of a type. In  Section \ref{sec:alg} we concentrate attention on types from the positions of one-sorted algebraic logic. Section \ref{sec:lg-type} deals with the ideas of universal logical geometry which give rise to $LG$-types and their geometric description. We finish the paper with the list of  problems appearing in the context of previous considerations. 

{\bf Acknowledgements.} E.~Aladova was  supported  by the Minerva
foundation through the Emmy Noether Research Institute, by the
Israel Science Foundation and ISF center of excellence 1691/10.
The support of these institutions is gratefully appreciated.
E.~Plotkin is thankful for the support of the Minerva foundation
through the Emmy Noether Research Institute.

\section{Structures of algebraic logic}\label{sec:al}

We consider algebra and logic with respect to a given variety of
algebras $\Theta$. This point of view (cf. \cite{P9})  implies
some differences with the original notions introduced by P.~Halmos
(\cite{Halmos} -- \cite{Halmos_4},  see \cite{Leblanc} for
non-homogenious polyadic algebras). For the sake of convenience,
in this section we provide the  reader with all necessary
definitions. It will be emphasized that the transition from pure
logic to logic in $\Theta$ is caused by many reasons, and we would
like to distinguish the needs of universal algebraic geometry
among them.

  Denote by $\Omega$ the signature of operations in algebras from $\Theta$. Let $W(X)$ denote the free in $\Theta$ algebra over a non-empty set of variables $X$. In the meantime we assume that each $X$ is a subset of some infinite set of variables $X^0$.

We shall recall the well-known definitions of the existential and
universal quantifiers which are considered as new operations on
Boolean algebras (see \cite{Halmos}).

Let $B$ be a Boolean algebra. The mapping $\exists : B \to B$ is
called an {\it existential} quantifier if
\begin{enumerate}
\item
  $\exists (0) =0$,

\item
  $ a \le \exists (a) $,

\item
 $\exists (a \wedge \exists b) = \exists a \wedge \exists b$.
\end{enumerate}

\noindent The {\it universal} quantifier $\forall : B \to B$ is
defined dually:
\begin{enumerate}
\item
  $\forall(1) =1$,

\item
  $a \ge \forall (a)$,

\item
 $\forall (a \vee \forall b) = \forall a \vee \forall b$.
\end{enumerate}

\noindent Here the numerals $0$ and $1$ are zero and unit of the
Boolean algebra $B$ and $a,b$ are arbitrary elements of $B$.
Symbol $=$ means coincidence of elements in Boolean algebra, i.e.,
$a\le b$ and $b\le a$ is written as $a=b$, $a,b\in B$. The
quantifiers $\exists$ and $\forall$ are coordinated by: $\neg
(\exists a)=\forall (\neg a)$, i.e., $ (\forall a)=\neg(\exists
(\neg a))$ .

A pair $(B,\exists)$, where $B$ is a Boolean algebra and $\exists$
is the existential quantifier, is a {\it monadic} algebra (see
\cite{Halmos}).


\begin{defin}\label{Def_Quantifier}
A Boolean algebra $B$ is a
 {\it quantifier $X$-algebra} if a quantifier $\exists
 x \colon B \to B$ is defined  for every variable $x \in X$,
 and
 $$
  \exists x \exists y = \exists y \exists x,
 $$
 for every $x, y \in  X$.
\end{defin}

\begin{remark} See also the definition of diagonal-free cylindric algebras of Tarski e.a. \cite{HMT}.
\end{remark}

\begin{remark}
According to \cite{Halmos}, \cite{P9} a  Boolean algebra $B$ is a
 {\it quantifier $X$-algebra} if a quantifier $\exists
 (Y)\colon B \to B$ is defined  for every subset $Y \subset X$,
 and

 1. $  \exists (\varnothing)=I_B$,  the identity function on $B$,

 2. $ \exists (X_1\bigcup X_2)=\exists (X_1) \exists (X_2),$ where $X_1$, $X_2$ are subsets in $X$.
 \noindent

 If we restrict ourselves with finite nontrivial subsets of $X$, then these two definitions coincide, because  condition 2) implies commutativity of quantifiers, and, conversely, one can define $\exists(Y)=\exists y_1\cdots \exists y_k$, where $Y=\{y_1,\ldots, y_k\}$.

\end{remark}

 We shall consider also {\it quantifier $W(X)$-algebras $B$ with equalities}.
An equality in a quantifier $W(X)$-algebra is  symmetric,
reflexive and transitive  (see Definition \ref{Def_ExtendedBA})
predicate  $\equiv:W(X)\times W(X)\to B$ which takes a pair $w, w'
\in W(X)$ to the constant in $B$ denoted by $w \equiv w'$, subject
to condition:



1)  $(w_1\equiv w'_1 \wedge\ldots \wedge w_n \equiv w'_n) \le (w_1
\ldots w_n \omega \equiv w'_1\ldots w'_n \omega) $ where $\omega$
is an $n$-ary operation in $\Omega$.

We can speak about quantifier $W(X)$-algebras, assuming that  the
free in $\Theta$ algebra $W(X)$ uniquely corresponds  to each set
$X$. Suppose that the logical signature  is extended by symbols of
nullary operations $w \equiv w'$, where $w,w' \in W(X)$. Then

\begin{defin}\label{Def_ExtendedBA}
We call a Boolean algebra $B$  a quantifier $W(X)$-algebra  with
equalities (or an  extended Boolean algebra  over the free in
$\Theta$ algebra $W(X)$),
 if

1. There are defined quantifiers $\exists x$ for all $x \in X$ in
$B$ with $\exists x \exists y = \exists y \exists x$ for all $x,y
\in X$.

2. To every pair 
   $w,w' \in W(X)$ it corresponds a constant (called an equality) in  $B$, denoted
by $w \equiv w'$. Here,

2.1. $w_1\equiv w'_1 \le w'_1\equiv w_1$.

2.2. $w \equiv w$ is the unit of the algebra $B$.

2.3. $w_1\equiv w_2 \wedge w_2\equiv w_3 \le w_1\equiv w_3$.

2.4. For every $n$-ary operation $\omega\in \Omega$, where
$\Omega$ is a signature of the variety $\Theta$, we have
$$
w_1\equiv w'_1 \wedge \ldots \wedge w_n \equiv w'_n\le w_1 \ldots
w_n \omega \equiv w'_1 \ldots w'_n \omega.
$$
\end{defin}

\begin{remark}  Under homomorphisms of extended Boolean algebras each constant  $w \equiv w'$ goes to another constant of the same kind. Endomorphisms  of Boolean algebras leave constants $w \equiv w'$ unchanged.
\end{remark}

\begin{remark}
Condition 2.4 means that
 for every homomorphism $\mu:W(X)\to H$, where $H\in \Theta$, there is a
coordination of $\mu$ with all operations from $\Omega$. In other
words equalities respect all operations on $W(X)$.
\end{remark}


\begin{defin}\label{Def_OneSortHalmosAlg_1}
An algebra $\mf L=\mf L(X) $ is a Halmos algebra (one-sorted Halmos algebra) over  $W(X)$, 
$X$ is infinite 
 if:
\begin{enumerate}
\item
 $\mf L$ is an extended Boolean algebra.

\item
 The action of the semigroup $End(W(X))$ is defined on $\mf L$, so that for each $s\in
 End(W(X))$ there is the map
  $s_{*}: \mf L \to \mf L$ which preserves the Boolean structure of $\mf L$.

\item The identities controlling the interaction of $s_{*}$ with
quantifiers are as follows:

\begin{enumerate}
\item[3.1.]
 $s_{1*} \exists x a = s_ {2*} \exists x
a, \ a \in \mf L$, if $s_1( y) = s_2(y)$ for every $y
 \in X, \; y \neq
x$.

  \item[3.2.]
   $s_{*}\exists x a = \exists (s(x)) (s_*a),$ $\ a \in \mf L$,
if $ s(x) = y $ and $y $ 
is a variable which does not belong 
to the support 
 of $s(x')$, 
 for every 
$x' \in X$, and $x' \neq x$.

This condition means that $y$  does not participate in the
shortest expression of the element $s(x')\in W(X)$ through the
elements of $X$.

 \end{enumerate}

\item
   The identities controlling the interaction of $s_{*}$ with equalities are as follows:
  \begin{enumerate}
  \item[4.1]
   $s_{*}(w\equiv w')=(s(w)\equiv s(w'))$.

  \item[4.2]
   $(s^{x}_{w})_{*}a \wedge (w\equiv w')\le (s^{x}_{w'})_{*}a$,
   where $a\in \mf L$, and $s^{x}_{w}\in End(W(X))$ is defined by $s^{x}_{w}( x)=w$,
   and $s^{x}_{w}(x')=x',$ for $x'\ne x$.

 \end{enumerate}

\end{enumerate}
\end{defin}

\begin{remark}\label{rm:inf}
The set $X$ in the definition \ref{Def_OneSortHalmosAlg_1} must be
infinite
because otherwise $End(W(X))$ does not act on $B$ (see \cite{P9},
Chapter 8, Section 2 for the details) in the case of free Halmos
algebras. In general this condition is superfluous since we
require the action of the semigroup $End(W(X))$  on the algebra
$\mf L$.

 For the definition of support see \cite{P9}, Chapter 9, Section 1. 
\end{remark}

\begin{remark}
Definition \ref{Def_OneSortHalmosAlg_1} introduces algebras which
are very close to polyadic algebras of Halmos  ( see
\cite{Halmos})  defined over a set of variables $X$. The main
difference between these classes comes from the desire to
specialize an  algebraization of first order logic  to an
arbitrary variety of algebras $\Theta$. This means that instead of
action of the semigroup of transformations $End(X)$ of the set
$X$, we consider the action of the bigger semigroup $End(W(X))$ as
the semigroup of Boolean endomorphisms. We also consider
equalities of the type $w\equiv w'$ instead of the ones $x\equiv
y$ for polyadic algebras.
\end{remark}

\begin{remark}
Axioms 3.1 and 3.2 which look messy, are grounded on major examples of Halmos algebras. In particular, we will see that Halmos algebras of the kind $Hal_\Theta(H)$ (see Example \ref{ex:AlgebraBool}) satisfy these identities. Since these algebras generate the whole variety of Halmos algebras, 
every Halmos algebra should satisfy these identities. If instead
of $Hal_\Theta(H)$ we consider the Halmos algebra of formulas
$\widetilde \Phi$ (see below), then the identity 3.1. corresponds
to the well-known fact that it is possible to replace a quantified
variable in a formula by another one. The identity 3.2. has a
similar explanation (see \cite{Halmos}).

\end{remark}

\begin{remark} In \cite{Halmos}, \cite{P9} an equality in Halmos algebras is defined as a reflexive binary predicate which satisfies conditions 4.1. and 4.2. Then, it can be checked \cite{P9}, that this predicate is automatically symmetric and transitive. 

\end{remark}

\begin{example}\label{ex:AlgebraBool}We give an example of Halmos algebra which plays a crucial role in further considerations.
\end{example}

Let $X$ be any set (finite or infinite), $H$ an algebra in $\Theta$. Consider the set 
$Hom(W(X),H)$ of all homomorphisms from $W(X)$ to $H$.   Let
$Bool(W(X),H)$ be the Boolean algebra of all subsets $A$ in
$Hom(W(X),H)$. Our aim is to make it an extended Boolean algebra.

Define, first,  quantifiers $\exists x,$  $x \in X$ on
$Bool(W(X),H)$. We set $\mu \in \exists x A$ if and only if there
exists $\nu \in A$ such that $\mu(y) = \nu(y)$ for every $y \in
X$, $y\neq x$. It can be checked  that $\exists x$ defined in such
a way is, indeed, an existential quantifier.

Let us consider equalities of the form $w\equiv w'$, where
$w,w'\in W(X)$. Define the corresponding elements of the algebra
$Bool(W(X),H)$ as follows
$$
Val^{X}_{H}(w\equiv w')=\{ \mu \mid \mu (w)=\mu(w') \}.
$$
The set $Val^{X}_{H}(w\equiv w')$ is considered as an equality in
the algebra $Bool(W(X),H)$.

Thus, the algebra $Bool(W(X),H)$ is equipped with the structure of
an  extended Boolean algebra (we omit verification of the
necessary axioms).

Let  $X$ now be an infinite set. Define the action of the
semigroup $End(W(X))$ in $Bool(W(X),H)$. Every homomorphism $s\in
End(W(X))$ gives rise to a Boolean homomorphism
$$
s_{*}: Bool(W(X),H)\to Bool(W(X),H),
$$
defined by the rule: for each $A\subset Hom(W(X),H)$ the point
$\mu$ belongs to $s_{*}A$ if $\mu s\in A$.

The signature of a Halmos algebra for $Bool(W(X),H)$  is now completed, one can check  
that all axioms are satisfied and thus,  $Bool(W(X),H)$ is a
Halmos algebra. Denote it by $Hal^X_\Theta(H)$.

\bigskip

Our next aim is to define multi-sorted Halmos algebras. There are
many reasons to do that. Some of them are related to potential
applications of algebraic logic in computer science, but some have
purely algebraic nature. For instance, we need multi-sorted
variant of Halmos algebras in order to work with finite
dimensional affine spaces and to construct geometry related to
first order calculus in arbitrary $\Theta$.


Every multi-sorted algebra $D$ can be written as 
$D=(D_i, i \in \Gamma)$, where $\Gamma$ is a set of sorts, which
can be infinite, and $D_i$ is a domain of the sort $i$. We can
regard domains $D_i$ as algebras from some variety (for
definitions see \cite{Kurosh_Lectures}, \cite{Malcev}).

Every operation $\omega$ in $D$ has a specific type $\tau = \tau
(\omega)$. This notion generalizes the notion of  the arity of an
operation. In the multi-sorted case an operation $\omega$ of the
type $\tau = (i_1, \ldots, i_n;j)$ operates as a mapping $\omega:
D_{i_1} \times \ldots \times D_{i_n} \to D_j$. Homomorphisms of
multi-sorted algebras act component-wise and have the form $\mu =
(\mu _i, i \in \Gamma) : D\to D'$, where $\mu_i : D_i \to D'_i$
are homomorphisms of algebras  and, besides that, every $\mu$ is
naturally correlated with the operations $\omega$.

Subalgebras, quotient algebras, and cartesian products of
multi-sorted algebras are defined in the usual way. Hence, one can
define  varieties of multi-sorted
algebras. 
In every such a variety there exist free algebras over
multi-sorted sets, determined by multi-sorted identities.

It is worth noting  that categories and multi-sorted algebras are
tightly connected  \cite{Higgins}, \cite{MacLane}. So, define,
first, Halmos categories. Let $\Theta^0$ be the category of free
algebras of the variety $\Theta$.

\begin{definition}
A category $\Upsilon$ is a Halmos category if: 
\begin{enumerate}
\item Every its object has the form $\Upsilon(X)$, where
$\Upsilon(X)$ is an extended Boolean algebra in $\Theta$ over
$W(X)$. \item Morphisms are of the form $s_* : \Upsilon(X) \to
\Upsilon(Y)$, where every $s: W(X) \to W(Y)$ is a homomorphism in
$\Theta^0$, $s_*$ is the homomorphism of Boolean algebras and the
correspondence: $W(X)\to \Upsilon(X)$ and $s \to s_*$ determines a
covariant functor $\Theta^0 \to \Upsilon$. \item The identities
controlling the interaction of morphisms
with quantifiers and equalities repeat the ones from Definition \ref {Def_OneSortHalmosAlg_1}, where the endomorphisms $s$ from $End(W(X))$ are replaced by homomorphisms $s: W(X)\to W(Y)$. 
\end{enumerate}
\end{definition}
Now we are able to define  multi-sorted Halmos algebras associated
with Halmos categories. Consider an arbitrary $W(X)$ in $\Theta$
and take the signature $L_X = \{ \vee, \wedge, \neg, \exists x, x
\in X, M_X \}$.  Here $M_X$ is the  set of all equalities $w\equiv
w'$, $w,w'\in W(X)$ over the algebra $W(X)$. We treat equalities
from $M_X$ as nullary operations.  We  add all $s=s^{XY}: W(X) \to
W(Y)$ to all $L_X$, where $X, Y\in \Gamma$, treating them as
symbols of unary operations (under unary we mean that these
operations of the type $(X,Y)$  use just one argument). Denote the
new signature by $L_\Theta$. So,
$$
L_\Theta= \{ \vee, \wedge, \neg, \exists x, x \in X, M_X, X\in
\Gamma,   s=s^{XY}  \}
$$
The signature $L_\Theta$ is a multi-sorted signature and consists
of all one-sorted signatures $L_X$, where $X$ runs $\Gamma$, and
of all $s$.

For the aims of logical geometry 
we assume that $\Gamma$ is the set of all finite subsets of the
infinite set $X^0$.

\begin{remark}\label{rm:inf} This condition on the domains $\Gamma$ is not necessary for the definition of Halsmos algebras  and made exclusively for geometric needs. Halmos algebras can be defined  for  various choice of domains. For example, the one-sorted Halmos algebra from Definition \ref{Def_OneSortHalmosAlg_1} corresponds to the signature $L_X = \{
\vee, \wedge, \neg, \exists x, x \in X, M_X , s \}$, where $s:
W(X) \to W(X)$, and $X=X^0$ is an infinite set.
\end{remark}


Consider further algebras $\Upsilon = (\Upsilon_X, X \in \Gamma)$.
Every $\Upsilon_X$ is an algebra in the
 signature $L_X$ and a unary operation (mapping) $s_*: \Upsilon_X \to \Upsilon_Y$  corresponds to every
$s: W(X) \to W(Y)$.

\begin{definition}\label{ha:ms}
We call an algebra $\Upsilon= (\Upsilon_X, X \in \Gamma)$ in the
signature $L_\Theta$ a Halmos algebra, if

1. Every $\Upsilon_X$ is an extended Boolean algebra in the
signature $L_X$.

2. Every mapping $s_*: \Upsilon_X \to \Upsilon_Y$ is 
  a homomorphism of Boolean
algebras.

3. The identities, controlling interaction of operations $s_*$
with quantifiers and equalities are the same as in the definition
of Halmos categories.

4. Let $s: W(X) \to W(Y)$, $s': W(Y) \to W(Z)$, and let $u\in
\Upsilon_X$. Then $s'_*(s_*(u))=(s's)_*(u)$.
\end{definition}

\noindent It is clear that each Halmos category $\Upsilon$ can be
viewed as a Halmos algebra and vice versa.

\begin{remark} The choice of  $\Theta$ gives rise to some conditions all $s_*$ have to satisfy.
\end{remark}

Now we shall  construct two major examples of multi-sorted Halmos
algebras. The first one mimics the construction of one-sorted
Halmos algebra from Example \ref{ex:AlgebraBool}.

{\bf 1.} Our aim is to define the Halmos category
$Hal_{\Theta}(H)$. Assume that we have a class of sets $X_i$,
$X_i\in \Gamma$. Objects of this category  are extended Boolean
algebras $Bool(W(X_i),H)$ from Example \ref{ex:AlgebraBool}, for
various $X_i\in \Gamma$. Morphisms
$$
s_{*}:Bool(W(X_i),H)\to Bool(W(X_j),H),
$$
 are defined as follows:
$$
\mu\in s_{*}A  \Leftrightarrow \widetilde s(\mu)=\mu s\in A,
$$
where $\mu:W(X_j)\to H$, $A\subset Hom(W(X_i),H)$, $s:W(X_i)\to
W(X_j)$, and $\widetilde s: Hom(W(X_j),H)\to Hom(W(X_i),H)$. Here
$\widetilde s$ is viewed as a morphism of the category of affine
spaces.  In other words, $s_{*}A=(\widetilde s)^{-1}A$. A morphism
$s_{*}$ is automatically a homomorphism of Boolean algebras. The maps $s_{*}$ 
 are correlated with quantifiers and equalities, see
\cite{Plotkin_AG} for details. Moreover, there is a covariant
functor: $ \Theta^{0}\to Hal_{\Theta}(H). $ Hence,
$Hal_{\Theta}(H)$ is a Halmos category.


The category $Hal_{\Theta}(H)$ gives rise to a multi-sorted
($\Gamma$-sorted) Halmos algebra, denoted by
$$
Hal_{\Theta}(H)=(Bool(W(X_i),H),X_i\in\Gamma).
$$
Each component here is the extended Boolean algebra. The
operations in $Hal_{\Theta}(H)$ are presented by the operations in
each component $Bool(W(X_i),H)$ and unary operations corresponding
to
 morphisms
$$
s_{*}: Bool(W(X_i),H)\to Bool(W(X_j),H).
$$

{\bf 2.}  Another important example of multi-sorted Halmos algebra
is presented by algebra $\widetilde \Phi=(\Phi(X), X\in \Gamma)$
of  first order formulas with equalities. It turns out that
geometrical aims forces to consider multi-sorted variant of
algebraization of first order calculus and consider multi-sorted,
in a special sense, formulas. We shall return to this discussion
at the end of the section.

Consider once again the signature $L_\Theta= \{ \vee, \wedge,
\neg, \exists x, x \in X, M_X, X\in \Gamma,   s=s^{XY}  \} $,
where  $M_X$ is the  set of all equalities $w\equiv w'$, $w,w'\in
W(X)$ over $W(X)$ treated as  nullary operations, and
$s=s^{XY}:W(X)\to W(Y)$ are symbols of  unary operations.

First, we construct the algebra  $\widetilde \Phi$ in an explicit
way. Denote by $M=(M_X, X \in \Gamma)$ the multi-sorted  set of
equalities with the components $M_X$.

 Each equality $w\equiv w'$ is a formula of the length
zero, and of the sort  $X$ if $w\equiv w'\in M_{X}$. Let $u$ be a
formula of the length $n$ and the sort $X$. Then the formulas
$\neg u$ and $\exists x u$ are the formulas of the same sort $X$
and  the length $(n+1)$. Further, for the given $s:W(X)\to W(Y)$
we have the formula $s_{*}u$ with the length $(n+1)$ and the sort
$Y$. Let now $u_1$ and $u_2$ be formulas of the same sort $X$ and
the length $n_1$ and $n_2$ accordingly. Then the formulas $u_1\vee
u_2$ and $u_1\wedge u_2$ have the length $(n_1+n_2+1)$ and the
sort $X$. In such a way, by induction, we define lengths and sorts
of arbitrary 
formulas.

Let $\mathfrak L^{0}_{X}$ be the set of all formulas of the sort
$X$. Each $\mathfrak L^{0}_{X}$ is an algebra in the signature
$L_{X}$ and 
$$
\mathfrak L^{0}=(\mathfrak L^{0}_{X}, X\in \Gamma)
$$
is an algebra in the signature $L_{\Theta}$. By construction,
algebra $\mathfrak L^{0}$ is the absolutely free algebra of
formulas over equalities (i.e. over nullary operations) concerned
with the variety of algebras $\Theta$.

Denote by  $ \tilde \pi$ the congruence in $\mathfrak L^{0}$
generated by the identities of Halmos algebras from Definition
\ref{ha:ms} (see also their list in Definition
\ref{Def_OneSortHalmosAlg_1}) and define the Halmos algebra of
formulas as
$$
\widetilde \Phi=\mathfrak L^{0}/ \tilde \pi;
$$
It can be  written as $\widetilde \Phi=(\Phi(X), X\in \Gamma)$,
where
$$
\Phi(X)=\mathfrak L^{0}_{X}/ \tilde \pi_{X},
$$
where each $\Phi(X)$ is an extended Boolean algebra of the sort
$X$ in the signature $L_X$. The algebra $\widetilde\Phi$ is,
obviously, the free algebra in the variety of all  multi-sorted
Halmos algebras associated with the variety of algebras $\Theta$,
with the set of free generators $M=(M_X, X\in \Gamma)$. Denote
this variety by $Hal_\Theta$.

\begin{remark}\label{rm:ms}
One can show \cite{Plotkin_IsotAlg}, that if we factor out 
component-wisely 
 the algebra $\mathfrak L^{0}$ by the 
many-sorted 
 Lindenbaum-Tarski congruence, then we get the same
algebra $\widetilde \Phi$. This observation provides a bridge
between syntactical and semantical description of the free
multi-sorted Halmos algebra.
\end{remark}

\begin{remark}\label{rm:ms1} 
To the contrary of the one-sorted case, the described construction does not give much practical information about multi-sorted formulas. Indeed, suppose we consider a one-sorted algebra $\Phi(X)$. 
 Let us
pick up an arbitrary element $u$ from  $\Phi(X)$. We can consider this element  as a mirror in the one-sorted Halmos algebra $\Phi(X)$ of a first order formula constructed on the base of the equality predicate. 
Looking at the element we can deduce the structure of the
corresponding formula.

The existence of  operations $s:W(X)\to W(Y)$ breaks this
intuition in many-sorted case. If an element $u$ has the sort $X$
and thus belong to $\Phi(X)$, then we cannot represent explicitly
the element
$s_{*}u$ from $\Phi(X)$ in terms of equalities, connectives, and quantifiers in $\Phi(X)$. This means that we cannot trace the structure of an arbitrary element from $\Phi(X)$.
\end{remark}

Fortunately, there exists a way out from the difficulty described
in Remark \ref{rm:ms1}. If we were to know what  the algebras
which constitute the variety $Hal_\Theta$ are, then we could
calculate the image of any element from $\Phi(X)$ in algebras from
$Hal_\Theta$. The following theorem yields that this is the case
in our situation.

\begin{theorem}[\cite{Plotkin_AG}]\label{th:tdf} The variety $Hal_\Theta$ of multi-sorted Halmos algebras is generated by all algebras $Hal_\Theta(H)$, where $H\in\Theta$.
\end{theorem}

Theorem  \ref{th:tdf}, in fact, gives us another definition for
the algebra $\widetilde \Phi$,
which can be considered as a free algebra in the variety generated
by  algebras $Hal_\Theta(H)$. 
 This allows us to study properties
of $\widetilde \Phi$ using the very concrete algebra $$
Hal_{\Theta}(H)=(Bool(W(X),H),X\in\Gamma)
$$
 as a model.
 Recall that we have defined the image of equalities from $M_{X}$ in $Bool(W(X),H)$ by:
 $$
Val^{X}_{H}(w\equiv w')=\{ \mu \mid \mu (w)=\mu(w') \}.
$$
This means that there is the map

$$Val_{H}: M\to Hal_\Theta(H).$$
Since equalities $M=(M_X, X\in\Gamma)$   freely generate  the free
multi-sorted Halmos algebra $\widetilde\Phi$, the map  $Val_{H}$
can be extended from generators to the homomorphism of
multi-sorted Halmos algebras
$$
 Val_{H}: \widetilde\Phi \to Hal_\Theta(H).
 $$
Since $\widetilde\Phi =(\Phi(X), X\in \Gamma)$, where each
component $\Phi(X)$ is an extended Boolean algebra, the
homomorphism   $
 Val_{H}$ induces homomorphisms
$$
 Val^{X}_{H}: \Phi(X) \to Bool(W(X),H),
 $$
of the one-sorted extended Boolean algebras. This allows
  us to calculate the value of each element from $\Phi(X)$ in  $Bool(W(X),H)$. Note that the values of elements of the form $s_\ast u$ are calculated as follows.
  Take $s:W(X)\to W(Y)$ and consider the formula $s_\ast u$, where $u\in \Phi(X)$. By definition, $s_\ast u$ belongs to  $\Phi(Y)$. Since $Val_H$ is a homomorphism, then

   $$Val^{Y}_H(s_\ast u)=s_\ast(Val^{X}_H u).$$

In the next sections we shall put all  this staff in the context
of affine spaces in arbitrary varieties. Replacing usual equations
by logical formulas we arrive at the field of logical geometry
which is
 much more complicated than the ordinary
equational geometry.

\section{Structures of universal algebraic geometry}\label{sec:galois}

Let us begin with the very classical setting (cf. \cite{Shafarevich}). Let $K$ be a field and 
$T=\{f_1,\ldots, f_m\}$ be a set polynomials   in the polynomial
algebra $K[X]=K[x_1,\ldots,x_n]$. Consider the affine space $K^n$
with points $\bar a= (a_1,\ldots, a_n)$, $a_i\in K$ and define the
Galois correspondence between ideals $T$ in $K[X]$ and algebraic
sets $A$ in $K^n$:
$$
T'_K=A=\{\bar a \ | \ f_i(\bar a)=0, \ \mbox{ for \ all }  f_i\in
T \},
$$
$$
A'_K=T=\{f_i\in K[X] \ | f_i(\bar a) =0, \ \mbox{ for \ all } \
\bar a\in A\},
$$
In this correspondence geometric objects: curves, surfaces,
general algebraic sets appear as zero loci of polynomials in the
algebra $K[X]$.

In order to generalize this situation to arbitrary varieties of
algebras, consider the variety $Com-K$ of commutative, associative
algebras with unit over the field $K$. Then the algebra $K[X]$ is
the free algebra in this variety and polynomials $f_i$ are just
elements of free algebra. Consider the field $K$ and its
extensions as algebras in this variety. Consider elements $\bar
a=(a_1,\ldots,a_n)$ of the affine space $K^n$ as functions $\bar
a:K[X]\to K$ defined by $\bar{a}(x_i)=a_i$, $i=1,\ldots, n$. Using
this vocabulary we can define the Galois correspondence and
geometric objects not in $Com-P$ but in arbitrary $\Theta$.

Let $\Theta$ be an arbitrary  variety and $H$ be an algebra in
$\Theta$. This algebra takes the role of the field $K$, hence the
affine space has to be of the form $H^n$. Let $W(X)$ be the free
algebra over $X$, $X=\{x_1,\ldots,x_n\}$. This is the place were
equations are situated and thus it plays the role of
$K[x_1,\ldots,x_n]$. The natural bijection $\alpha: Hom(W(X),H)\to
H^n$ allows us to consider the set of homomorphisms $Hom(W(X),H)$
as the affine space and its elements as the points of the affine
space. Let the  point $\mu\in Hom(W(X),H)$ be induced by a map
$\mu: X\to H$. Then it corresponds the point $\bar a=(a_1, \ldots,
a_n)$ in $H^n$, where $a_i=\mu(x_i)$. This correspondence gives
rise to kernels of points $\mu$ of the affine space. We define the
kernel $Ker(\mu)$ of the point $\mu$   as the kernel of the
homomorphism $\mu:W(X)\to H$.

Let  $T$ be a system of equations of the form $w\equiv w'$,
$w,w'\in W(X)$ which we treat as  a system of formulas of the form
$w\equiv w'$ on $W(X)$. Since
 $w$ and $w'$ are formulas in $W(X)$, then
$w=w(x_1, \ldots, x_n)$, $w'=w'(x_1, \ldots, x_n)$.

\begin{definition}
A point $\bar a=(a_1, \ldots, a_n)\in H^n$ is a solution of
$w\equiv w'$ in the algebra $H$ if
 $w(a_1, \ldots, a_n)=w'(a_1, \ldots, a_n)$.
A point $\mu\in Hom (W(X),H)$ is a solution of $w\equiv w'$ if
 $w(\mu(x_1), \ldots, \mu(x_n)) = w'(\mu(x_1), \ldots, \mu(x_n))$.
 \end{definition}

   The equality
 $w(\mu(x_1), \ldots, \mu(x_n)) = w'(\mu(x_1), \ldots, \mu(x_n))$
means that the pair $(w, w')$ belongs to $Ker(\mu)$. In other words, a point $\mu$ is a solution of the equation  $w\equiv w'$ if this formula belongs to the kernel of the point $\mu$. 
 Thus we say that  $w\equiv w'$ belongs to the kernel of a point if and only if the pair $(w, w')$ belongs to this kernel.
 The kernel
$Ker(\mu)$ is a congruence of the algebra $W(X)$, and the quotient
algebra $W(X) / Ker(\mu)$ is defined.

Let now $T$ be a system of equations in $W(X)$ 
and $A$ a set of points in $Hom(W(X),H)$. 
 Set the Galois
correspondence by
$$
T'_H = A = \{\mu : W(X) \to H \ | \ T \subset Ker(\mu)\}
$$
$$
A'_{H}=T=\{(w\equiv w')  \ |\ (w, w')\in \bigcap_{\mu\in A}
Ker(\mu)\}.
$$

\begin{definition} A set $A$ in the affine space $Hom(W(X),H)$ is called an algebraic set if there exists a system of
equations $T$ in $W(X)$ such that each point $\mu$ of $A$
satisfies all equations  from $T$. A congruence $T$ in $W(X)$ is
called $H$-closed if there exists $A$ such that $ A'_{H}=T$.
\end{definition}
We can rewrite the Galois correspondence through the values of
formulas:

$$
T'_H = A=\bigcap_{(w,w') \in T} Val^X_H (w\equiv w').
$$
$$
A'_H=T =  \{w\equiv w' \ |\ A \subset Val^X_H (w\equiv w')\}.
$$

The geometry obtained via this correspondence is an equational
geometry grounded on algebra $H$ in $\Theta$. However, there are
no reasons to restrict ourselves with equational predicates
looking at the images of the formulas in the affine space. We can
look at arbitrary first order  formulas as at equations, and since
arbitrary formulas are the elements of $\widetilde\Phi=(\Phi(X),
X\in \Gamma)$,  we shall replace in all consideration the free
algebra $W(X)$ by the extended Boolean algebra $\Phi(X)$.

The sets of equations are defined as arbitrary subsets in
$\Phi(X)$, the finite dimensional affine space $H^n$ is the same
as in equational case, and it remains to define the geometric
objects, that is the images of the formulas $u\in \Phi(X)$ in the
Galois correspondence. This can be done because, as we know,  the
equalities $M_X$, $X\in\Gamma$ represent the free generators of
$\widetilde\Phi$ and, thus the value homomorphism $Val^X_H$ can be
extended from equalities to arbitrary formulas $u\in \Phi(X)$.

Let $\mu: W(X)\to H$ be a point. Along with the classical kernel
$Ker(\mu)$ we define its logical kernel.

\begin{defin}
A formula $u\in \Phi(X)$ belongs to the logical kernel $LKer(\mu)$
of a point $\mu$ if and only if $\mu\in Val^{X}_{H}(u)$.
\end{defin}

It can be verified that the logical kernel $LKer(\mu)$ is always a
Boolean ultrafilter of $\Phi(X)$ \cite{Plotkin_IsotAlg}.

Since we consider each formula $u\in \Phi(X)$ as an "equation" and
$Val^{X}_{H}(u)$ as a value of the formula $u$ in the algebra
$Bool(W(X),H)$, then  $Val^{X}_{H}(u)$ is a set of points
$\mu:W(X)\to H$ satisfying the "equation" $u$. We call
$Val^{X}_{H}(u)$  solutions of the equation $u$. We also say  that
the  formula $u$ holds true in the
algebra $H$ at the point $\mu$. 

 {\it We call the obtained geometry associated to an arbitrary variety $\Theta$ and $H\in\Theta$ the logical geometry}.

 In order to
establish in this case the Galois correspondence   we shall
replace the kernel $Ker(\mu)$ by the logical kernel
$LKer(\mu).$ Let $T$ be a set of formulas in $\Phi(X)$ 
and $A$ a set of elements in $Bool(W(X),H)$. 
 Define
$$
T^{L}_{H}=A=\{ \mu :W(X)\to H \mid T\subset LKer(\mu) \},
$$
$$
A^{L}_{H}=T=\bigcap_{\mu\in A} LKer(\mu)
$$
The same Galois correspondence can be rewritten as
$$
T^{L}_{H}=A=\bigcap_{u\in T} Val^{X}_{H}(u).
$$
$$
A^{L}_{H}=T=\{ u\in\Phi(X) \mid A\subset Val^{X}_{H}(u) \}.
$$

\begin{definition} A set $A$ in the affine space $Hom(W(X),H)$ is called an elementary set if there exists a system of
formulas $T$ in $\Phi(X)$ such that each point $\mu$ of $A$
satisfies all formulas  from $T$. In other words,
$A=A^{LL}_{H}=T^{L}_{H}$ is fulfilled for elementary sets.
\end{definition}

\begin{definition}
A set of formulas $T\subset \Phi(X)$ such that
$T=T^{LL}_{H}=A^{L}_{H}$ is called {\it an $H$-closed Boolean
filter} in $\Phi(X)$.
\end{definition}

\begin{remark} The set of formulas $T$ which defines an elementary set $A$ can be infinite.
\end{remark}

\begin{remark} Elementary sets in the model theory are usually  called 
definable sets. 
 Since in the geometrical approach they are
tightly connected with elementary theories, we use the term
"elementary set" instead of "definable set".
\end{remark}

\begin{remark}\label{rem:diff} The formulas from $T\subset \Phi(X)$ may contain free generators  from different  $X_i$, $i\in \Gamma$. For example, the formula
$$
u=s_\ast(y_1\equiv y_2)\vee (x_3\equiv x_4),
$$
where $X=\{x_1,x_2,x_3,x_4\}$, $Y=\{y_1,y_2\}$ and $s(y_1)=x_1$,
$s(y_2)=x_2$, belongs to $\Phi(X)$.
\end{remark}

\section{Model theoretic types}\label{sec:mod}

In this  section we have to recall, first, the well-known
definitions from model theory. In our exposition, we mainly follow
the standard model theory course by \cite{Marker}, see also
\cite{ManinZilber}, \cite{Poizat}, etc. We assume that the precise
definition of an $\mathbb L$-structure is known. Basically, an
$\mathbb L$-structure is a pair $(\mathbb L, M)$, where $\mathbb
L$ is a language and $M$ is a set, called the domain of the
structure. Any language may contain functional symbols, symbols of
relations, and special symbols called constants. Given an $\mathbb
L$-structure, all these symbols are  interpreted (realized) on the
domain $M$. So any  $\mathbb L$-structure can be considered as a
triple $(\mathbb L, M, f)$, where $f$ is an interpretation
function.

 Formulas of $\mathbb L$ are  built inductively from atomic formulas, using the symbols of $\mathbb L$,
symbols of variables $x_1, x_2, \dots$, the equality symbol
$\equiv$, the Boolean connectives $\wedge$, $\vee$, $\neg$, the
quantifiers $\exists$ and $\forall$, and parentheses $(\ $, $)$.
We suppose that the interpretation of symbol $\equiv$ is always
equality on $M$.

A variable $x$ {\it occurs freely} in a formula $u$ if it is not
bounded by quantifiers
$\exists x$ or $\forall x$. A formula $u$ is called {\it a sentence} (or a 
closed 
formula) if it has no free variables. 
If $u(x_1, \dots , x_n)$ is a formula in free variables $x_1,
\dots , x_n$ then its closure $\bar u$ is any sentence produced
from $u$ by bounding all free variables by quantifiers. 


Let $\mathbb M$ be an $\mathbb L$-structure. For an $\mathbb
L$-formula $u$ one writes $\mathbb M \models u$ to say that the
value of $u$ under the interpretation $f$ is true.
  The value 
  ("true" or "false") under interpretation $f(x_i)=a_i$, $i=1,\ldots,m$, $a_i\in M$ 
  of a formula $u=u(x_1,\ldots,x_m)$ is defined inductively, using Tarski schema. Each $\mathbb L$-sentence is either true or
false on 
the whole 
$\mathbb M$. Let $u(x_1, \dots , x_n)$ be a formula in free
variables $x_1, \dots , x_n$  which means that all occurrences of
other variables in this formula are bounded. If $u(x_1, \dots ,
x_n)$ is  a formula with free variables $x_1, \dots , x_n$ and
$\bar a=(a_1,\dots
,a_n)\in M^{n}$ then we write $\mathbb M \models u(a_1, \dots ,a_n)$ for the true formula under interpretation $f(x_i)=a_i$. In this case we say that $u$ is satisfiable on $\mathbb M$. 

\begin{defin}\label{Def_Theory_4}
A set $T$ of $\mathbb L$-sentences is called an $\mathbb
L$-theory. $\mathbb M$
 is a model of the theory $T$ if $\mathbb M\models u$ for all  
  $u\in
T$. A theory is satisfiable if it has a model.
\end{defin}

Suppose that $\mathbb M$ is an $\mathbb L$-structure and
$A\subseteq M$. Let $\mathbb L_{A}$ be the language obtained by
adding to $\mathbb L$ constant symbols for each $a\in A$. We can
naturally view $\mathbb M$ as an $\mathbb L_A$-structure by
interpreting the new symbols in the obvious way. Let
$Th_{A}(\mathbb M)$ be the set of all $\mathbb L_{A}$-sentences
true in $\mathbb M$, that is the $\mathbb L_{A}$-{\it theory} of
the model $\mathbb M$.

\begin{defin}\label{Def_Theory_2}
If $\mathbb L$ is a first order language, then $Th_{A}(\mathbb M)$
is called the {\it elementary theory} of $M$.
\end{defin}

\begin{defin}\label{Def_MT-Type}
Let $P=\{u_i(x_1,\ldots,x_n)\}$  be a set of $\mathbb
L_{A}$-formulas in free variables $x_1,\dots , x_n$. We call $P$
an $n$-type  (partial $n$-type) if $P\cup Th_{A}(\mathbb M)$ is
satisfiable. We say that $P$  is a complete $n$-type if $u\in P$
or $\neg u \in P$ for all $\mathbb L_{A}$-formulas $u$ with free
variables from $x_1, \dots ,x_n$.
\end{defin}

So, the data for a type $P$ is a structure $\mathbb M$ and a subset of constants $A\subseteq M$. 
If $\mathbb M$ is any $\mathbb L$-structure, $A\subseteq M$, and
$\ol a= (a_1,\dots ,a_n)\in M^{n}$, let $tp^{\mathbb M}(\ol
a/A)=\{ u(x_1,\dots ,x_n)\in\mathbb L_{A}: \mathbb M\models
u(a_1,\dots ,a_n) \}$. Then, $tp^{\mathbb M}(\ol a/A)$ is a
complete $n$-type.


\begin{defin}\label{Def_MT-Type_1}
We say that a complete n-type $P$ is realized in $\mathbb M$ if
there is $\ol a= (a_1,\dots ,a_n)\in M^{n}$ such that
$P=tp^{\mathbb M}(\ol a/A)$.
\end{defin}

Denote the sets of all complete realizable $n$-types over $M$  by
$S^n_A(\mathbb M)$. In case $A=M$ we denote this set by
$S^n(\mathbb M)$.
\begin{problem}\label{prob:1}
Suppose that for two structures  ${\mathbb M_1}$ and ${\mathbb
M_2}$ the sets of complete realizable $n$-types $S^n(\mathbb M_1)$
and $S^n(\mathbb M_2)$ coincide for every $n$. What can be said
about ${\mathbb M_1}$ and ${\mathbb M_2}$?  How far are these
structures  from being isomorphic?
\end{problem}
\begin{remark}
Topologically, this question is very close to the following one:
suppose two structures have isomorphic Stone spaces (i.e., the
spaces of  complete realizable $n$-types $S^n(\mathbb M)$) for
each $n$. What can be said about relations between the structures
in this case?
\end{remark}


Problem \ref{prob:1} is a generalization of the problem about elementary equivalence of structures. Loosely speaking we ask how distant can  algebraic structures be 
if not only their logical descriptions coincide, but coincide also
the logical descriptions of particular elements from these
structures. This question can be specialized to specific varieties
of algebras $\Theta$ and to specific algebras in $\Theta$.

\section{Algebraization of model theoretic types}\label{sec:alg}

Define an algebraization of the notion of type.   Let $X^0$ be an
infinite set of variables. Let $H$ be an algebra from a variety of
algebras $\Theta$. Let the set of  constants equal $H$, that is we
consider  algebras $G$ from the variety $\Theta^H$ of
$H$-algebras. For example, if $\Theta$ is the variety of
commutative and associative rings with the unit and $K$ is a
field, then $\Theta^K$ is the variety of algebras over the field
$K$.

In our case, the free algebras in $\Theta^H$  have the form
$W(X^0)=W'(X^0)\ast H$, where $W'(X^0)$ is the free algebra in
$\Theta$ and $\ast$ stands for the free product in $\Theta$.

Let $\Phi(X^0)$ be the one-sorted Halmos algebra of formulas
associated with the variety $\Theta^H$. Recall that $\Phi(X^0)$ is
constructed in the following way.  We consider the signature
consisting  of symbols of Boolean connectives, existential
quantifiers $\exists x$, $x\in X^0$, equalities of the form
$w(x_1,\ldots,x_n)\equiv w'(x_1,\ldots,x_n)$, where $w, w'$
belongs to $W(X)$,  $X$ runs all finite subsets of $X^0$,  and
symbols of operations $s: W(X)\to W(X)$, for every $X$. Let us
take the absolutely free algebra over equalities in this
signature. The quotient of this algebra by the Lindenbaum-Tarski
congruence is $\Phi(X^0)$. The pair $(\Phi(X^0), H)$ plays the
role of $\mathbb L_M$-structure   ${\mathbb M}$, where $M=H$.

Now we recall the Galois correspondence from the previous section
in the case when $\widetilde\Phi=(\Phi(X), X\in \Gamma)$ is a
one-sorted Halmos algebra $\Phi(X^0)$, $X^0$ is infinite. Let $T$
be a set of formulas in $\Phi(X^0)$. We have
$$
T^{L}_{H}=A=\{ \mu :W(X)\to H \mid T\subset LKer(\mu) \},
$$
$$
A^{L}_{H}=T=\bigcap_{\mu\in A} LKer(\mu)
$$
In particular, $u\in T$ if and only if $A\subset Val^{X}_{H}(u)$.

 Let $X=\{x_1,\ldots,x_n\}$ be a finite subset  in $X^0$. We shall define $X$-$MT$-type ($MT$-type for short) of the point $\mu\in Hom(W(X),H)\cong H^n$.

 For each point $\mu: W(X)\to H$ consider the set of points $A_\mu$ defined by: a point  $\nu:W(X^0)\to H$ belongs to $A_\mu$ if $\nu(x)=\mu(x)$ for $x\in X$ and $\nu(y)$ is an arbitrary element in $H$. Define
$$
T_\mu=(A_\mu)^L_{ H}=\bigcap_{\nu\in A_\mu}LKer(\nu).
$$
In other words $ T_\mu$  is the set of all formulas $u\in
\Phi(X^0)$ which hold on the points from $A_\mu$. This means that
$u\in T_\mu$ if $A_\mu\subset Val^{X^0}_H(u)$. Since every logical
kernel is an ultrafilter, the set $T_\mu$ is a filter.

\begin{definition}\label{def-MT} We call the filter $
T_\mu$ an $MT$-type of the point $\mu$.
\end{definition}

\begin{remark}\label{rem:type}
Let us compare Definitions \ref{Def_Theory_4} --
\ref{Def_MT-Type_1} and Definition \ref{def-MT}. In the definition
\ref{def-MT} we consider an $MT$-type of the point $\bar
a=(a_1,\ldots, a_n)$, where $\mu(x_i)=a_i$, $a_i\in H$ for $x_i\in
X$,  as the set of all formulas $u$ which hold true on the point
$\mu$ (i.e., on the point $\bar a$). Therefore, the type of a
point in our definition is always a filter.

On the other hand, by the definition \ref{Def_MT-Type_1} the type
of the point $tp^{\mathbb H}(\bar a)=tp^{\mathbb H}( \mu) $, where
$\mu (x_i)=a_i$, $i=1,\ldots,n$ is the set  of the satisfiable in
the point $\mu$  formulas of the form
$u=u(x_1,\ldots,x_n,y_1,\ldots, y_k)$, where only $x_i$ are free
variables. This is a subset of $T_\mu$ and thus an $MT$-type
$T_\mu$ is somewhat bigger than the corresponding $tp^{\mathbb H}(
\mu)$.


\end{remark}

\begin{remark}\label{rem:etheory}
The similar situation holds with the definition of the elementary
theory of an algebra $H$.

We will consider elementary theory of $H$ as the set of all
formulas $u$ true in every point $\mu: Hom(W(X),H)$.

 On the other side, according to the modal-theoretic  Definition \ref{Def_Theory_4} the elementary theory of $H$ is smaller and consists of   closed formulas true in $H$. Since every formula $u$ true in $H$ is equivalent to its closure $\bar u$, then by abuse of language we use the same notation $Th(H)$ for the elementary theory of $H$ in both cases. So,
 $$
 Th(H)=\bigcap_{\mu} T_\mu,
 $$
\noindent where $\mu\in Hom(W(X),H)$.

This situation is typical for algebraic logic and geometry where the free variables do not play the same role as in logic and model theory. 
\end{remark}

Denote  the system of all $MT$-types $T_\mu$ of the algebra $H$ by
$S^X_H$. Here, $\mu:W(X)\to H$, and $X$ runs all finite subsets of
$X^0$.

Given finite subset $X\subset X^0$ and a point $\mu:W(X)\to H$,  define $s=s^\mu: W(X^0)\to W(X^0),$ where $W(X^0)=W'(X^0)\ast H$,  by letting $s (x_i)=\mu(x_i)$, if $x_i\in X$, 
and $s (y)=y$ for $y\in Y^0=X^0\setminus X$. Let $s^\mu_\ast:
\Phi(X^0)\to \Phi(X^0)$ be the corresponding map of Halmos
algebras. 


\begin{prop}\label{pr:1} A formula $u\in\Phi(X^0)$ belongs to $T_\mu$ if and only if $s^\mu_\ast u$ belongs to the elementary theory $Th(H)$.
\end{prop}

{\it Proof.} Let  $s^\mu_\ast u$ belong to the elementary theory
$Th(H)$. We shall prove that $u\in T_\mu$. Thus, we shall check
that  $A_\mu\subset Val^{X^0}_H(u)$. Let $\nu \in A_\mu$. Let
$\delta:W(X^0)\to H$ be an arbitrary point in $Hom (W(X^0), H)$.
Then, for $x_i\in X $, we have $\delta
s^\mu(x_i)=\delta(\mu(x_i))=\mu(x_i)$ since $\delta$ fixes
constants. Correspondingly,  $\delta s^\mu(y_i)=\delta (y_i)$.
Thus we can choose $\delta$ such that
$\delta s^\mu=\nu$ 
for any $\nu \in A_\mu$. 
 Since $s^\mu_\ast u \in Th(H)$, then
$\delta$ lies in $Val^{X^0}_H (s^\mu_\ast u)=s^\mu_\ast
Val^{X^0}_H ( u)$. The latter equality means,  by definition, that
$\delta s^\mu$ lies in $Val^{X^0}_H(u)$. Hence, $A_\mu\subset
Val^{X^0}_H(u)$.

Conversely,  let $u\in T_\mu$. We shall prove that $s^\mu_\ast u$
belongs to the elementary theory $Th(H)$. So we have to check that
any point $\delta$ satisfies $s^\mu_\ast u$. Consider $\delta
s^\mu$. This point belongs to $A_\mu$. Hence $\delta s^\mu$ lies
in $Val^{X^0}_H(u)$. This means that $\delta$ lies in $s^\mu_\ast
Val^{X^0}_H(u)=Val^{X^0}_H(s^\mu_\ast u)$. Thus, an arbitrary
point $\delta$ belongs to $Val^{X^0}_H(s^\mu_\ast u)$ and
$s^\mu_\ast u$ lies in $Th(H)$. $\qquad \square$


Let  $u=u(x_1,\ldots,x_n, y_1,\ldots, y_k)$ be a formula in
$\Phi(X^0)$ such that $x_i\in X$, $y_i\in Y$, and all occurrences
of $x_i$ are free, all occurrences of $y_i$ are bounded.
 We call such a formula special.  

Let $u$ be a special formula. It can be seen that $s^\mu_\ast u$
replaces all occurrences of free variables $x_i$ by the their
images $h_i\in H$ under the homomorphism $s^\mu$. Hence
$s^\mu_\ast u$ has all variables bounded, i.e., $s^\mu_\ast u$ is
a sentence.

   Any $MT$-type is complete with respect to special formulas. Indeed, let $u$ be a special formula and let $u\notin T_\mu$. Consider $\neg u$. We have $s^\mu_\ast (\neg u)=\neg s^\mu_\ast ( u)$. By Proposition \ref{pr:1}, $s^\mu_\ast ( u)$ does not hold in $H$. Since $s^\mu_\ast u$ is a sentence, the formula $\neg s^\mu_\ast ( u)$ holds in $H$. Hence, $s^\mu_\ast (\neg u)$ holds in $H$ and thus belongs to $Th(H)$. Then $\neg u\in T_\mu$ according to Proposition \ref{pr:1}.

   Suppose now that for two algebras $H_1$ and $H_2$ the sets $S^X_{H_1}$ and $S^X_{H_2}$ of   $MT$-types $T_\mu$ coincide.
   Every $MT$-type contains the corresponding model theoretic $n$-type, where $n=|X|$. So the problem \ref{prob:1} can be restated  as {\t what can be said about the closeness of algebras $H_1$ and $H_2$} if $T_\mu$ and $T_\nu$ coincide?

   From now on, one can build  the type theory from the positions of one-sorted algebraic logic. In the next section we consider a more geometric approach, related to multi-sorted logic and multi-sorted Halmos algebras.

 \section{Logically-geometric types}\label{sec:lg-type}

Let us take the free multi-sorted Halmos algebra of formulas
$\widetilde \Phi=(\Phi(X), X\in \Gamma)$, where all $X$ are
finite. Recall the necessary facts from the previous sections.

 There is the  value homomorphism of  multi-sorted Halmos algebras $ Val_{H}: \widetilde\Phi \to Hal_\Theta(H)  $, which induces homomorphisms  of extended Boolean algebras $
 Val^{X}_{H}: \Phi(X) \to Bool(W(X),H) $, where $
Hal_{\Theta}(H)=(Bool(W(X),H),X\in\Gamma) $. We can write $
Val_{H}=(Val^X_H, X\in \Gamma)$. For every $X$, the homomorphism
$Val^X_H$ gives rise to a major Galois correspondence of logical
geometry between $H$-closed congruences in $ \Phi(X)$ and
elementary sets in finite dimensional affine spaces $Hom(W(X),H):$

$$
T^{L}_{H}=A=\{ \mu :W(X)\to H \mid T\subset LKer(\mu) \},
$$
$$
A^{L}_{H}=T=\bigcap_{\mu\in A} LKer(\mu).
$$

Let $Th(H)=(Th^X(H), X\in \Gamma)$ be the multi-sorted
representation of the elementary theory of $H$. We call its
component $Th^X(H)$ the {\it $X$-theory of the algebra $H$}. We
have:
$$
 Ker(Val_H)=Th (H),
 $$
 $$
 Ker(Val^X_H)=Th^X(H).
 $$
 The key diagram which relates logic of different sorts in multi-sorted case is as follows:

 $$
\CD
\Phi(X) @> s_\ast >> \Phi(Y)\\
@V  \kp^X_H  VV @VV \kp^Y_H V\\
Bool(W(X),H) @>s_\ast>> Bool(W(Y),H)
\endCD
$$

Here the upper arrow represent the syntactical transitions in the
category $Hal_\Theta$, the lower level does the same with the
respect to semantics in $Hal_\Theta$, and the correlation is
provided by the value homomorphism.

Recall that a formula $u\in \Phi(X)$ belongs to the logical kernel
$LKer(\mu)$ of a point $\mu$ if and only if $\mu\in
Val^{X}_{H}(u)$, that is $u$ lies in $LKer(\mu)$ if a point $\mu$
satisfies the "equation" $u$. This is the Boolean ultrafilter,
which contains $Th^X(H)$. Indeed, if $u\in Th^X(H)$ then
$Val^H_X(u)=Hom(W(X),H)$. In particular, $\mu\in Val^X_H(u)$ and
$u\in LKer(\mu)$. Thus $Th^X(H)\subset LKer(\mu).$ Moreover,
$$
Th^X(H)=\bigcap_{\mu} LKer \mu.
$$

Define now the concept of an $LG$-type.

\begin{defin}
Every ultrafilter $T$ in the algebra $\Phi(X)$  containing
$Th^X(H)$  is called $X$-$LG$-type.
\end{defin}

\begin{defin}
A type $T$ is called $X$-$LG$-type of the algebra $H$, if there is
a point $\mu : W(X)\to H$ such that $T=LKer(\mu)$.
\end{defin}

In the latter case we also  say that the type $T$ is realized in
$H$. Since the elementary $X$-theory is contained in each
$LKer(\mu)$ then the elementary $X$-theory $Th^{X}(H)$ is
contained in each $X$-$LG$-type of $H$. Denote  the system of all
$X$-$LG$-types of the algebra $H$ by $S^X(H)$.

Now we want to explore the geometrical nature of the Galois
correspondence. In algebraic geometry, the category of all
algebraic sets is an important invariant of the the algebra $H$.
In most cases, this category  is dual to the category of
coordinated algebras. We want to use similar ideas in the case of
logical geometry. The logical kernels take the role played by the
radical ideals in classical geometry and the roles of  closed
congruences in the universal one. So, the types of the points
represented by the logical kernels may have  similar impact to
logical geometry and may be involved in the similar
algebraically-geometric ideas.

Two algebras $H_1$ and $H_2$ are called {\it geometrically
equivalent} ($AG$-equivalent for short) (see
\cite{Plotkin_IsotAlg}, \cite{Plotkin_7Lec}) if for every finite
$X$ and $T$ in $W(X)$ we have
$$
T_{H_1}^{''}= T_{H_2}^{''}.
$$

\begin{defin} [\cite{PZ}]
Algebras $H_1$ and $H_2$ are called logically equivalent
($LG$-equivalent for short) if for every finite $X$ and $T$ in
$\Phi(X)$ we have
$$
T^{LL}_{H_1}= T^{LL}_{H_2}.
$$
\end{defin}

It can be seen (see \cite{PZ}), that if two algebras $H_1$ and
$H_2$ are logically equivalent then they are {\it elementary
equivalent} (i.e., $Th(H_1)=Th(H_2)$). The converse statement is
not true.

\begin{defin}[\cite{PZ}]
Algebras $H_1$ and $H_2$ in $\Theta$ are called $LG$-isotyped, if
for any finite $X$, every $X$-$LG$-type of the algebra $H_1$ is an
$X$-$LG$-type of the algebra $H_2$ and vice versa.
\end{defin}

Thus, the algebras $H_1$ and $H_2$ are $LG$-isotyped if
$S^X(H_1)=S^X(H_2)$ for every $X\in \Gamma$. This coincidence
clearly implies that they are elementary equivalent.

So, we have the geometric notion of logical equivalence of
algebras which generalizes geometric equivalence, and the model
theoretic notion of $LG$-isotypeness. Both of them imply
elementary equivalence.  The following theorem shows that these
two notions coincide.

\begin{theorem}[\cite{Plotkin_IsotAlg}]
Algebras $H_1$ and $H_2$ are $LG$-equivalent if and only if they
are $LG$-isotyped.
\end{theorem}

One can define the category of algebraic sets $K_\Theta(H)$ and
the category of elementary sets $LK_\Theta(H)$. The objects of
$K_\Theta(H)$ are of the form $(X,A)$, where $A$ is an algebraic
set in $Hom(W(X),H)$. If we take for $A$ the elementary sets, then
we are getting to the category of elementary sets $LK_\Theta(H)$.
The morphisms are of the form
$$
[s] : (X,A) \to (Y,B).
$$
Here $s:W(Y) \to W(X)$ is a morphism in the category $\Theta^0$.
The corresponding $\tilde s: Hom(W(X),H) \to Hom(W(Y),H)$ should
be coordinated with $A$ and $B$ by the condition: if $\nu \in
A\subset Hom(W(X),H)$, then $\tilde s (\nu) \in B\subset
Hom(W(Y),H)$. Then the induced mapping $[s]:A \to B$ we consider
as a morphism $(X,A) \to (Y,B)$.

The category $K_\Theta(H)$ is a full subcategory in
$LK_\Theta(H)$. It is known that if two algebras $H_1$ and $H_2$
are geometrically equivalent, then the categories of algebraic
sets $K_\Theta(H_1)$ and $K_\Theta(H_2)$ are isomorphic. A similar
fact is valid with respect to categories of elementary sets.
Namely,

\begin{theorem}[\cite{PZ}]
If the algebras $H_1$ and $H_2$ are $LG$-isotyped then the
categories $LK_\Theta(H_1)$ and $LK_\Theta(H_2)$ are isomorphic.
\end{theorem}

\section{Problems}\label{probl}

In Sections \ref{sec:alg} and \ref{sec:lg-type} we described
$MT$-types and $LG$-types. Now we want to compare these notions.

Recall that $MT$-types are defined for points $\mu: W(X)\to H$ of
the affine space $Hom(W(X),H)$. However, the formulas from any
$MT$-type $T_\mu$ lie in the algebra of formulas $\Phi(X^0)$,
where $X^0$ is an infinite set. It is important to note, that the
algebra $H$ from the given variety of algebras $\Theta$ is treated
as the algebra of constants.

In the case of $LG$-types, we consider finite sets $X$ 
 in $X^0$ and the multi-sorted algebra of formulas $\widetilde \Phi=(\Phi(X), X\in \Gamma)$, where all $X$ are finite.
 The $X$-$LG$-type of the point $\mu: W(X)\to H$ is $LKer(\mu)$, which is calculated in the algebra $\Phi(X)$. This is one of the differences in two approaches. We shall also remember  that  the formulas from $T\subset \Phi(X)$ may contain free generators  from different  $X$, where  $X\in \Gamma$ (see Remark \ref{rem:diff}).

\begin{prob}
Compare $MT$-isotypeness and $LG$-isotypeness. In other words, are
there algebras $H_1$ and $H_2$ such that they are $MT$-isotyped
but not $LG$-isotyped, or such that they are $LG$-isotyped but not
$MT$-isotyped?
\end{prob}

Problems \ref{pr:free} and \ref{pr:free1} are devoted to
$LG$-types.

\begin{prob}\label{pr:free}
Let $F_n$ be a free group of the rank $n>1$ and $H$ be a finitely
generated group. Is it true that if $F_n$ and $H$ are
$LG$-isotyped then they are isomorphic?
\end{prob}

\begin{prob}\label{pr:free1}
Are there $LG$-isotyped groups $H_1$ and $H_2$ such that $H_1$  is
finitely generated and $H_2$ is an arbitrary non finitely
generated group?
\end{prob}

C.~Perin and R.~Sklinos \cite{PerinSklinos} proved that if for a
non-abelian free group there is the equality $T_{\mu}=T_{\nu}$
then $\mu=\sigma\nu$ for some automorphism $\sigma$
of $H$. 

\begin{prob}
What are the varieties $\Theta$ such that for arbitrary free
algebra $H=W(X)$ from $\Theta$ the equality $T_\mu=T_\nu$ implies
$\mu=\sigma\nu$?
\end{prob}

Similar question for $LG$-types and free groups is of great
interest.

\begin{prob}\label{prob:22}
Is it true that for a given free non-abelian  group  the equality
$LKer(\mu)=LKer(\nu)$ implies $\mu=\sigma\nu$?
\end{prob}

Problem \ref{prob:22} has positive solution for the case of free
abelian groups (G.~Zhitomirski, unpublished).

Note that the group of automorphisms of an algebra $H$ acts on the
affine space $Hom(W(X),H)$, and each elementary set is invariant
under this action. If for the algebra $H$ there are only a finite
number of $Aut(H)$-orbits in $Hom(W(X),H)$ for every $X$, then
there are only finite number of realizable $LG$-types  in
$\Phi(X)$.
 It can be shown that for  free abelian groups of the exponent $p$ this property is satisfied.  It would be  interesting to look for non-abelian examples.

\begin{prob}
Find examples of algebras $H$ such that for every $X$ there are
only a finite number of $Aut(H)$-orbits in $Hom(W(X),H)$.
\end{prob}

\end{document}